\documentclass[11pt]{article}
\usepackage{amsmath}
\usepackage{mathrsfs}
\usepackage{amsfonts}

\usepackage{amsfonts, amsmath, amssymb}
\usepackage{amssymb,amsfonts,amsmath,
latexsym, epsfig,cite, psfrag,eepic,colordvi,color}
\usepackage{amscd,graphics}
\usepackage{lineno}

\newcommand{\qed}{\hfill $\Box $}
\newcommand{\pf}{\noindent {\bf Proof.} }
\newcommand{\rem}{\noindent {\bf Remark} }

\newtheorem{theorem}{Theorem}[section]
\newtheorem{lemma}[theorem]{Lemma}
\newtheorem{coro}[theorem]{Corollary}

\newtheorem{problem}[theorem]{Problem}

\begin{document}

\title{An Extension of Cui-Kano's Characterization Problem on Graph Factors
}

\author{ Hongliang Lu\thanks{Corresponding email: luhongliang215@sina.com (H. Lu)}
\\ {\small  School of Mathematics and Statistics}
\\ {\small Xi'an Jiaotong University, Xi'an 710049, PR China}
}

\date{}

\maketitle



\begin{abstract}
Let $G$ be a graph with vertex set $V(G)$ and let $H:V(G)\rightarrow
2^N$ be a set function associating with $G$.
 An $H$-factor of graph $G$ is a
spanning subgraphs~$F$ such that
$$d_F(v)\in H(v)\hspace{4em}\hbox{for every }v\in V(G).$$
Let $f:V(G)\rightarrow N$ be an even integer-valued function such
that $f\geq 4$ and let $H_f(v)=\{1,3,\ldots,f(v)-1, f(v)\}$ for
$v\in V(G)$. In this paper, we   investigate $H_f$-factors of graphs
$G$ by using Lov\'asz's structural descriptions.  Let $o(G)$ denote
the number of odd components of $G$.
 We show that if one
of the following conditions holds, then $G$ contains an
$H_f$-factor.
\begin{itemize}
\item[$(i)$] $o(G-S)\leq f(S)$ for all $S\subseteq V(G)$;

\item[$(ii)$]   $|V(G)|$ is odd, $d_G(v)\geq f(v)-1$ for all $v\in V(G)$
and $o(G-S)\leq f(S)$ for all
 $\emptyset\neq S\subseteq V(G)$.
\end{itemize}
As a corollary, we show that if a graph  $G$ with odd order and
minimum degree $2n-1$ satisfies
$$o(G-S)\leq 2n|S|\hspace{4em}  \mbox{ for all }
 \emptyset\neq S\subseteq V(G),$$
then $G$ contains an $H_n$-factor. 
In particular, we make progress on the characterization problem for
a special family of graphs proposed by Akiyama and Kano.
\end{abstract}

\section{Introduction}

  All graphs in this paper are simple. Let $G$ be a graph
with vertex set $V(G)$ and edge set $E(G)$. We denote the degree of
$v$ in $G$ by $d_G(v)$. The minimum degree in graph $G$ will be
denoted by $\delta(G)$ and the maximum degree by $\triangle(G)$.
The subgraph induced by the set $S$ is denoted by $G[S]$. The number
of components of graph $G$ is denoted by $\omega(G)$ and the number
of odd components of~$G$  by $o(G)$.  Let $E_G(S,T)$ denote the set
of edges of graph $G$   with one end in $S$ and the other end in $T$
and $e_G(S,T)=|E_G(S,T)|$.
 The join $G = G_1 + G_2$, is the graph obtained from two
vertex disjoint graphs $G_1$ and $G_2$  by joining each vertex in
$G_1$ to every vertex in $G_2$.

Let~$H$~be a function associating a subset of $\mathbb{Z}$ to each
vertex of~$G$. A spanning subgraph $F$ of graph $G$ is called   an
$H$-{\it factor} of $G$ if
\begin{align}\label{cond_fac}
d_{F}(x)\in H(x) \hspace{4em} \mbox{ for every vertex $x\in V(G)$}.
\end{align}
By specifying $H(x)$ to be an interval or a special set, an
$H$-factor becomes an $f$-factor, an $[a, b]$-factor or a $(g,
f)$-factor, respectively.

Let~$F$~be a spanning subgraph of~$G$. Following
Lov\'asz~\cite{Lovasz}, one may measure the ``deviation'' of~$F$
from the condition~\eqref{cond_fac} by
\begin{equation}\label{eq_deviation}
\nabla_H(F)=\sum_{v\in V(G)} \min\bigl\{|d_F(v)-h|\,\colon\,h\in
H(v)\bigr\}.
\end{equation}
Moreover, the ``solvability'' of~\eqref{cond_fac} can be
characterized by
\[
\nabla(H)=\min\{\nabla_H(F)\,\colon\,\mbox{$F$ is a spanning
subgraph of~$G$}\}.
\]
The subgraph $F$ is said to be {\em $H$-optimal} if
$\nabla_H(F)=\nabla(H)$. It is clear that $F$ is an $H$-factor if
and only if $\nabla_H(F)=0$, and any $H$-factor (if exists) is
$H$-optimal. Let
\[
Q=\{h_1,h_2,\ldots,h_m\},
\]
where $h_1<h_2<\cdots<h_m$. Then $Q$ is called an \emph{allowed set}
if each of the gaps of $Q$ has at most one integer, i.e.,
\[
h_{i+1}-h_i\le2\hspace{3em}\hbox{for all }1\le i\le m-1.
\]
A set function $H$  associating with $G$ is called an \emph{allowed
set function} (following~\cite{Lovasz}) if $H(v)$ is an allowed set
for all $v\in V(G)$.
%

 Lov\'asz \cite{Lovasz} showed that if $H$ is not an allowed
set, then the decision problem of determining whether a graph has an
$H$-factor is known to be $NP$-complete. Cornu\'ejols~\cite{Cor88}
provided the first polynomial algorithm for the problem with~$H$
allowed.

A special case of $H$-factor problem is the so-called
\emph{$(1,h)$-odd factor   problem}, i.e., the problem with
\[
H(v)=\{1,\,3,\ldots,\,h(v)-2,\,h(v)\},
\]
where $h:V(G)\rightarrow N$ be an odd function. For a constant odd
integer $n\geq 1$, if $h(x)=n$ for all $x\in V(G)$, then $(1,h)$-odd
factor is called \emph{$(1,n)$-odd factor}.
 The first investigation of the $(1,n)$-odd factor problem
is due to Amahashi~\cite{Ama85}, who gave a Tutte type
characterization for graphs having a global odd factor.

\begin{theorem}[Amahashi]\label{thm_Amahashi}
Let $n$ be an odd integer. A graph $G$ has an $(1,n)$-odd factor if
and only if
\begin{equation}\label{cond_Amahashi}
o(G-S)\le n\,|S|\hspace{4em}\hbox{for all subsets } S\subset V(G).
\end{equation}
\end{theorem}
For general odd value functions~$h$, Cui and Kano~\cite{CK88}
established a Tutte type theorem.
\begin{theorem}[Cui and Kano, \cite{CK88}]\label{CK88}
Let $h:V(G)\rightarrow N$ be odd value function. A graph $G$ has an
$(1,h)$-odd factor if and only if
\begin{equation}\label{cond_Amahashi}
o(G-S)\le h(S)\,\hspace{4em}\hbox{for all subsets } S\subset V(G).
\end{equation}
\end{theorem}

 Noticing the form of the condition~\eqref{cond_Amahashi},
they asked the question of characterizing graphs~$G$ in terms of
graph factors such that
\begin{equation}\label{cond_CK}
o(G-S)\le 2n\,|S|\hspace{4em}\hbox{for all subsets }S\subset V(G).
\end{equation}
Motivated by Cui-Kano's problem, Lu and Wang \cite{LuWang} consider
the degree prescribed subgraph problem for the special prescription
\begin{equation}\label{Hn}
H_n=\{1,3,\ldots,2n-1,2n\}.
\end{equation}

\begin{theorem}[Lu and Wang,\cite{LuWang}]\label{luwang}
Let $G$ be a connected graph. If
\begin{align}
o(G-S)\leq 2n\,|S|\hspace{4em}\hbox{for all subsets }S\subset V(G),
\end{align}
then $G$ contains an $H_n$-factor.

\end{theorem}
The condition of Theorem \ref{luwang} implies that $|V(G)|$ is even.
Let $H_n^*=H_n\cup\{-1\}$. For odd order graph, they obtained the
following result (for convenience, the definition of
$H_n^*$-critical graph will be introduced in Section
\ref{H-critical}).
\begin{theorem}[Lu and Wang,\cite{LuWang}]\label{thm_OddOrder}
Let $G$ be a connected graph of odd order. Suppose that
\begin{equation}\label{cond_nonempty}
o(G-S)\le 2n|S|\hspace{3em}\hbox{for all }\ \ \emptyset\ne S\subset
V(G).
\end{equation}
Then either $G$ contains an $H_n$-factor, or $G$ is
$H_n^*$-critical.
\end{theorem}

The condition~\eqref{cond_Amahashi}   implies that
 the graph is even order. For odd order graph, Akiyama and Kano
  propose the following problem (see also ~\cite[Problem~(6.14)]{AK11} ).

\begin{problem}[Akiyama and Kano, \cite{AK11}]\label{Aki_Ka}
Let $G$ be a connected graph and $h:V(G)\rightarrow N$ be an even
integer-valued function. If $G$ satisfies
\begin{align}
o(G-S)\leq h(S) \hspace{3em}\hbox{for all }\ \emptyset\ne S\subset
V(G),
\end{align}
what factor or property does $G$ has?
\end{problem}

 Let $f\geq 4$ be an even integer-value function and
let $H_f:V(G)\rightarrow 2^N$ be an set function such that
$H_f(v)=\{1,3,\ldots,f(v)-1,f(v)\}$ for $v\in V(G)$.
  Motivated by
Akiyama-Kano's problem, we   investigate the structure of graphs
without $H_f$-factor by using Lov\'asz's $H$-factor structure theory
\cite{Lovasz}. We obtain the following result, which is an extension
of Theorem \ref{luwang}.

\begin{theorem}\label{main_even}
Let $G$ be a graph with even order. If
\begin{align}\label{main_even_eq}
o(G-S)\leq f(S) \hspace{3em}\hbox{for all } S\subset V(G),
\end{align}
then $G$ contains an $H_f$-factor.
\end{theorem}

The inequality (\ref{main_even_eq}) also implies that $|V(G)|$ is
even. For odd order graph, we solve Problem \ref{Aki_Ka} and obtain
a stronger result than Theorem \ref{thm_OddOrder}.
\begin{theorem}\label{main_odd}
Let $G$ be a  connected graph with odd order. Suppose that
$d_G(v)\geq f(v)-1$ for all $v\in V(G)$. If
\begin{align}
o(G-S)\leq f(S) \hspace{3em}\hbox{for all }\ \emptyset\ne S\subset
V(G),
\end{align}
then $G$ contains an $H_f$-factor.
\end{theorem}

\begin{coro}\label{coro}
Let $n\geq 2$ be an integer and let $G$ be a connected graph with
odd order and minimum degree $2n-1$. If
\begin{align}\label{Coro_eq1}
o(G-S)\leq 2n|S| \hspace{3em}\hbox{for all }\ \emptyset\ne S\subset
V(G),
\end{align}
then $G$ contains an $H_n$-factor.
\end{coro}

\rem~\textbf{1:}  In Corollary \ref{coro},   the conditions
``$\delta(G)\geq 2n-1$" is sharp. Let $K_{2n-1}$ denote the complete
graph of order $2n-1$. Take $2n-2$ disjoint copies of $K_{2n-1}$.
Add a new vertices $v$ and connect  two vertices in each copy of
$K_{2n-1}$ to the new vertex $v$. This results a  connected graph
$G$ with odd order $(2n-2)(2n-1)+1$ and minimum degree $2n-2$. It is
easy to show that
\begin{align}
o(G-S)\leq 2n|S| \hspace{3em}\hbox{for all }\ \emptyset\ne S\subset
V(G).
\end{align}
Now we show that $G$ contains no $H_n$-factor. Otherwise, suppose
that $G$ contains an $H_n$-factor $F$. By parity, $K_{2n-1}$
contains no $H_n$-factors and so  $F$ contains exactly an edge from
a copy of $K_{2n-1}$. Then we have $d_F(v)=2n-1\notin H_n$, a
contradiction.

\rem~\textbf{2:} In Corollary \ref{coro}, the condition
(\ref{Coro_eq1}) is not necessary for the existence of an
$H_n$-factor in a graph. Let $m\geq 2n+2$ be an even integer.
Consider the graph
$$G=K_1+m K_{2n+1}$$
obtained by linking a vertex  $v$ to all vertices in $2n + 1$ copies
of the complete graph $K_{2n+1}$. Clearly, $G$ is a graph with odd
order and minimum degree $2n+1$. It is  easy to verify that $G$
contains an $H_n$-factor. However, taking the subset $S$ to be the
single vertex $v$, we see that the condition (\ref{Coro_eq1}) does
not hold for G.



\section{On $H$-critical Graphs}\label{H-critical}

In this section, we study $H$-factors of graphs based on Lov\'asz's
structural description to the degree prescribed subgraph problem.
 Denote by $I_H(v)$ the set of vertex
degrees in all $H$-optimal subgraphs of graph $G$, i.e.,
\[
I_H (v)=\{d_F(v)\,\colon\,\hbox{all $H$-optimal subgraphs $F$}\}.
\]
Comparing the set $I_H(v)$ with $H$, one may partition the vertex
set $V(G)$ into four classes:
\begin{align*}
C_H&=\{v\in V(G)\,\colon\,I_H(v)\subseteq H(v)\},\\[5pt]
A_H&=\{v\in V(G)-C_H\,\colon\,\min I_H(v)\ge \max H(v)\},\\[5pt]
B_H&=\{v\in V(G)- C_H\,\colon\,\max I_H(v)\le \min H(v)\},\\[5pt]
D_H&=V(G)- A_H- B_H- C_H.
\end{align*}
It is clear that the $4$-tuple $(A_H,B_H,C_H,D_H)$ is a pairwise
disjoint partition of $V(G)$. We call it the {\em $H$-decomposition}
of~$G$. In fact, the four subsets can be distinguished according to
the contributions of their members to the
deviation~\eqref{eq_deviation}. A graph $G$ is said to be {\em
$H$-critical} if it is connected and $D_H=V(G)$.  For
non-consecutive allowed set function, the only necessary condition
of $H$-critical graph is given by Lov\'asz \cite{Lovasz}. In this
paper, we obtain a sufficient condition for $H$-critical graph.

We write $MH(x)=\max H(x)$ and $mH(x)=\min H(x)$ for $x\in V(G)$.
For $S\subseteq V(G)$, let $MH(S)=\sum_{x\in S}MH(x)$ and
$mH(S)=\sum_{x\in S}mH(x)$.  By the definition of $A_H,
B_H,C_H,D_H$, the following holds:
\begin{itemize}
\item[$($I$)$] for every $x\in B_H$, there exists an $H$-optimal graph $F$ such
that $d_F(x)<mH(x)$;

\item[$($II$)$] for every $x\in A_H$, there exists an $H$-optimal graph $F$ such
that $d_F(x)>MH(x)$;

\item[$($III$)$] for every $x\in D_H$, there exists an $H$-optimal graph $F$ such
that $d_F(x)<MH(x)$ and other $H$-optimal graph $F'$ such that
$d_F(x)>mH(x)$.

\end{itemize}

  Lov\'asz \cite{Lovasz} gave the following
properties.
\begin{lemma}[Lov\'asz,\cite{Lovasz}]\label{lem_interval}
If $G$ is a simple graph, then $I_H(v)$ is an interval for all $v\in
D_H$.
\end{lemma}

\begin{lemma}[Lov\'asz,\cite{Lovasz}]\label{lem_gap}
The intersection $I_H(v)\cap H(v)$ contains no consecutive integers
for any vertex $v\in D_H$.
\end{lemma}

Given an integer set $P$ and an integer $a$, we write $P-a=\{p-a\ |\
p\in P\}$. Let $C$ be a connected induced subgraph of $G$ and
$T\subseteq V(G)-V(C)$. Let $H_{C,T}:V(C)\rightarrow 2^N$ be a set
function such that $H_{C,T}(x)=H(x)-e_G(x,T)$ for all $x\in V(C)$.
\begin{lemma}[Lov\'asz,\cite{Lovasz}]\label{lem_D_critical}
Every component $R$ of $G[D_H]$ is $H_{R,B_H}$-critical and if $F$
is $H$-optimal, then $F[V(R)]$ is $H_{R,B_H}$-optimal.
\end{lemma}

\begin{lemma}[Lov\'asz,\cite{Lovasz}]\label{lem_def}
If $G$ is $H$-critical, then $\nabla(H)=1$.
\end{lemma}

\begin{lemma}[Lov\'asz,\cite{Lovasz}]\label{lem_necessary}
For any $H$-optimal graph $F$, $E_G(B_H,B_H\cup C_H)\subseteq E(F)$,
and $E_G(A_H,C_H\cup A_H)\cap E(F)=\emptyset$.
\end{lemma}

\begin{theorem}[Lov\'asz,\cite{Lovasz}]\label{lem_defici}\label{thm1}
$\nabla(H)=\omega(G[D_H])+\sum_{v\in
B_H}(mH(v)-d_{G-A_H}(v))-\sum_{v\in A_H}MH(v)$.
\end{theorem}

In the proof of main theorems, we need the following two technical
lemmas.

\begin{lemma}\label{lem_misses}
Let $F$ be an $H$-optimal subgraph. For every component $R$ of
$G[D_H]$, $F$ misses at most an edge of $E_G(V(R),B_H)$.
\end{lemma}

\pf Let $F$ be an $H$-optimal subgraph of $G$. We write
$\tau_H=\omega(G[D_H])$ and $G[D_H]=C_{1}\cup\cdots\cup
C_{\tau_{H}}$. Since $C_i$ is $H_{C_i,B_H}$-critical, then $C_{i}$
contains no $H_{C_i,B_H}$-factors. So if  $def_{F}(C_{i})=0$,  then
$F$ either misses at least an edge of $E(C_{i},B_H)$ or contains at
least an edges of $E(C_{i},A_H)$. Let $\tau_{B}$ denote the number
of components of $G[D_H]$ such that $F$ misses at least an edge of
$E(C_{i},B_H)$ and $\tau_{A}$ denote the number of the components of
$G[D]$  such that $F$ contains at least an edge of $E(C_{i},A_H)$.
Let $\tau_c$ denote the number of components $C_i$ of $G[D_H]$ such
that $F$ contains at least one edge of $E(C_{i},A_H)$ and misses at
least one edge $E(C_{i},B_H)$. Then we have
\begin{align*}
\nabla_{H}(F)&\geq \tau_{H}-\tau_{A}-\tau_{B}+\tau_c+\sum_{x\in
A_H\cup
B_H}\min \{|r- d_{F}(x)|\ |\ r\in H(x)\}\\
&\geq \tau_{H}-\tau_{A}-\tau_{B}+\tau_c+\sum_{x\in
A_H}(d_{F}(x)-MH(x))+\sum_{x\in
B_H}(mH(x)-d_{F}(x))\\
&\geq
\tau_{H}-\tau_{A}-\tau_{B}+\tau_c+(e_{F}(A_H,B_H)+\tau_{A}-MH(A_H))+\sum_{x\in
B_H}(mH(x)-d_{F}(x))\\
&=\tau_{H}-\tau_{B}+\tau_c+(e_{F}(A_H,B_H)-MH(B_H))+\sum_{x\in
B_H}(mH(x)-d_{F}(x))\\
&\geq
\tau_{H}-\tau_{B}+\tau_c+(e_{F}(A_H,B_H)-MH(A_H))+(mH(B_H)-(e_{F}(A_H,B_H)+\sum_{x\in
B_H}d_{G-A_H}(x)-\tau_{B}))\\
&=\tau_{H}(A_H,B_H)+\tau_c+mH(B_H)-MH(A_H)-\sum_{x\in
B_H}d_{G-A_H}(x)\geq \nabla(H).
\end{align*}
Since $\nabla_{H}(F)=\nabla(H)$, then  we obtain $\tau_c=0$ and
\begin{align*}
\sum_{x\in B_H}d_{F}(x)=e_{F}(A_H,B_H)+\sum_{x\in
B_H}d_{G-A_H}(x)-\tau_{B_H},
\end{align*}
which implies that $F$ misses at most an edge from $C_i$ to $B$.
This completes the proof for $1\leq i\leq \tau_H$.  \qed

\begin{lemma}\label{non_consecu}  Let $G$ be a graph and let $H:V(G)$ be an allowed set function. If
 $MH(v)-1\in H(v)$ and $d_G(v)\geq MH(v)-1$ for all $v\in V(G)$, then $G$ is not
$H$-critical.
\end{lemma}

\pf By contradiction, we firstly assume that $G$ is  $H$-critical.
Let $F$ be an $H$-optimal subgraph of $G$ such that $E(F)$ is
 maximal.

 Since $G$ is $H$-critical
 and $F$ is $H$-optimal, then by Lemma \ref{lem_def}, we have
$d_F(v)\leq MH(v)+1$ for all $v\in V(G)$. We claim that there exists
a vertex $x\in V(G)$ such that $d_F(x)= MH(x)+1$. Otherwise, suppose
that $d_F(v)\leq MH(v)$ for all $v\in V(G)$. Then there exists a
vertex $v\in V(G)$ such that $d_{F}(v)\notin H(v)$ and so $0\leq
d_{F}(v)\leq f(v)-2$. Hence there exists an edge  $e \in E(G)-E(F)$,
which is
 incident with vertex $v$. Then $F\cup \{e\}$ is also
$H$-optimal, contradicting to the maximality of $F$. Thus there
exists a vertex $x\in V(G)$ such that $d_F(x)=MH(x)+1$. Since
$I_{H}(x)$ is an interval and $I_{H}(x)\cap H(x)$ contains no two
consecutive integers, then we have $\min I_{H}(x)\geq MH(x)$,
contradicting to $x\in D_{H}$.

This completes the proof. \qed

\begin{coro}
Let $n\geq 2$ be an integer and let $G$ be a graph. If
$\delta(G)\geq 2n-1$, then $G$ is not $H_n$-critical.

\end{coro}

\section{The Proof of  Theorems \ref{main_even} and \ref{main_odd}}

In this section, we assume that $f:V(G)\rightarrow Z^+$ be an even
integer-valued function such that $f \geq 4$ and
$H_f(v)=\{1,3,\ldots,f(v)-1,f(v)\}$ for all $v\in V(G)$.  

\begin{theorem}\label{lem_main}
Let $G$ be a graph and let $A_{H_f}$, $B_{H_f}$, $C_{H_f}$ and
$D_{H_f}$ be defined as above. Then
\begin{itemize}
\item[$(a)$] $E_G(B_{H_f},B_{H_f}\cup C_{H_f})=\emptyset$;

\item[$(b)$]
For every component $R$ of $G[D_{H_f}]$,
$|V(R)|+|E_G(V(R),B_{H_f})|\equiv1\pmod 2$;

\item[$(c)$] every component $R$ of $G[D_{H_f}\cup B_{H_f}]$ is odd.

\end{itemize}

\end{theorem}

\pf Firstly, we prove (a) by contradiction. Suppose that there
exists an edge $e\in E_G(B_{H_f},B_{H_f}\cup C_{H_f})$. Without loss
of generality, we assume that $e=uv$ and $u\in B_{H_f}$.  For any
$H$-optimal graph $F$, by Lemma \ref{lem_necessary}, $e\in E(F)$ and
so $d_F(v)\geq 1$, contradicting to the definition of $B_{H_f}$.
This completes the proof of (a).

Secondly, we prove (b). By Lemma \ref{lem_D_critical}, $R$ is
$H_{R,B_{H_f}}$-critical. For simplicity, we write
$H_R=H_{R,B_{H_f}}$. We claim that $f(u)-e_G(u, B_{H_f})\notin
I_{H_R}(u)$ for all $u\in V(R)$. Otherwise, suppose that there
exists a vertex $x\in V(R)$ such that $f(x)-e_G(x, B_{H_f})\in
I_{H_R}(x)$.
 By Lemma \ref{lem_gap}, $I_{H_R}(x)\cap H_R(x)$
does not contain two consecutive integers and so we have
$f(x)-1-e_G(u, B_{H_f})\notin I_{H_R}(x)$. By Lemma
\ref{lem_interval},   $I_{H_R}(x)$ is an interval and so we have
$\min I_{H_R}(x)\geq MH_R(x)$, contradicting to the definition of
$H$-critical graphs.
Hence $I_{H_R}(u)\subseteq [0,f(u)-1-e_G(u, B_{H_f})]$.
 Let $F$ be an $H_f$-optimal graph and $F^*=F[V(R)]$.
 By Lemma \ref{lem_D_critical},
 $F^*$ is an $H_R$-optimal subgraph of graph $R$.
 Furthermore, by Lemma \ref{lem_D_critical}, $R$ is $H_R$-critical and so there
exists a vertex $x\in V(R)$ such that $d_{F^*}(x)\notin H_R(x)$ and
$d_{F^*}(y)\in H_R(y)$ for all $y\in V(R)-x$.

Hence for every vertex $y\in V(R)-x$,  $d_{F^*}(y)\equiv
f(y)-1-e_G(y,B)\pmod 2$ and $d_{F^*}(x)\equiv
f(x)-2-e_G(x,B_{H_f})\pmod 2$. Then
\begin{align*}
\sum_{v\in V(R)}d_{F^*}(v)&\equiv
\sum_{y\in V(R)-x}(f(y)-1-e_G(y,B_{H_f}))+f(x)-2-e_G(x,B_{H_f})\pmod 2\\
&\equiv \sum_{y\in V(R)}e_G(y,B_{H_f})+|V(R)|-1,
\end{align*}
which implies
\begin{align*}
\sum_{y\in V(R)}e_G(y,B_{H_f})+|V(R)|\equiv 1\pmod 2.
\end{align*}
  This completes the proof of (b).

Finally, we prove (c).  We write
$B_{H_f}=\{v_1,\ldots,v_{|B_{H_f}|}\}$ and
$G[D_{H_f}]=C_1\cup\cdots\cup C_{\tau}$. For $1\leq i\leq |B_{H_f}|$
and $1\leq j\leq \tau$, we claim $e_G(v_i,V(C_j))\leq 1$. Otherwise,
suppose that there exists $v\in B_{H_f}$ and a component $C_{i}$ of
$G[D_{H_f}]$ such that $e_G(v,V(C_{i}))\geq 2$. For arbitrary
$H_f$-optimal graph $F$, by Lemma \ref{lem_misses}, then we have
$d_F(v)\geq 1$, contradicting $v\in B_{H_f}$.

Let $R$ be an arbitrary
 connected component of $G[D_{H_f}]$. Without loss of generality, we write
$V(R)=C_{1}\cup \cdots\cup C_{k}\cup B_1$, where
$B_1=\{y_1,\ldots,y_r\}$. Now we construct a   graph $R^*$ obtained
from $R$ by contracting $C_{i}$ to a   vertex $x_i$ for $1\leq i\leq
k$. By (a), $R^*$ is a bipartite graph.

\medskip
\textbf{ Claim 1.~} $R^*$ is a tree.
\medskip

Since $R$ is connected, then $R^*$ is connected. Now we show that
$R^*$ contains no cycles. Conversely, suppose that   $R^*$ contains
a cycle $x_1y_1,\ldots,x_my_m x_1$. We write
$W=C_{i_1}\cup\cdots\cup C_{i_m}$ and $B_2=\{y_1,\ldots,y_m\}$. By
Lemma \ref{lem_misses}, for any $H$-optimal graph $F$, $F$ contains
at least $m$ edges from $W$ to $B_2$. Now we claim that $d_F(v)= 1$
for all $v\in B_2$, otherwise, there exists a vertex $v\in B_2$ such
that $d_F(v)\geq 2$ contradicting to $v\in B_2\subseteq B_{H_f}$.
Since $F$ is an arbitrary $H_f$-optimal graph and $d_F(v)=1$ for all
$v\in B_2$, then we have $B_2\subseteq C_{H_f}$, a contradiction
again. This completes Claim 1.

Let $W^*=V(C_1)\cup \cdots\cup V(C_k)$. By Claim 1, $R^*$ is a tree,
which implies that $e_G(W^*,B_1)=k+r-1$. By (b),  we have
\begin{align*}
k&\equiv\sum_{i=1}^k (|V(C_i)|+e_G(V(C_i),B_{H_f}))\pmod 2\\
&=\sum_{i=1}^k
|V(C_i)|+e_G(W^*,B_1)\\
&=\sum_{i=1}^k |V(C_i)|+k+r-1,
\end{align*}
which implies
\begin{align*}
\sum_{i=1}^k |V(C_i)|+r\equiv 1\pmod 2.
\end{align*}
Hence $|V(R)|$ is  odd. This completes   the proof. \qed

By Theorems \ref{lem_defici} and \ref{lem_main}, we obtain the
following result.
\begin{coro}
If  graph $G$ contains no $H_f$-factors, then there exists  two
disjoint subsets $S,T$ of $V(G)$ such that
\begin{align*}
f(S)-|T|+\sum_{x\in T}d_{G-S}(x)-q(S,T)< 0,
\end{align*}
where $q(S,T)$ denote the number of components $C$ of $G-S-T$ such
that $|V(C)|+e_G(V(C),T)\equiv 1\pmod 2$.
\end{coro}

\noindent \textbf{Proof of Theorem \ref{main_even}.} Since $|V(G)|$
is even, by Theorem \ref{lem_main} (b), $G$ is not $H_f$-critical.
By Lemma \ref{lem_main} (a) and (c),
$E_G(B_{H_f},C_{H_f})=\emptyset$ and every component of
$G[D_{H_f}\cup B_{H_f}]$ is an odd component. We write
$\omega(G[D_{H_f}\cup B_{H_f}])=k$ and $G[D_{H_f}\cup
B_{H_f}]=R_1\cup \ldots\cup R_k$. Without loss of generality,
suppose that $V(R_i)=V(C_{i1})\cup \cdots\cup V(C_{i{r_i}})\cup B_i$
for $1\leq i\leq k$, where $C_{ij}$ is a component of $G[D_{H_f}]$
for $1\leq j\leq r_i$ and $B_i\subseteq B_{H_f}$. By Theorem
\ref{lem_defici}, Theorem \ref{lem_main} (c) and Claim 1 of Theorem
\ref{lem_main},
\begin{align*}
0<\nabla(H_f)&=\omega(G[D_{H_f}])+|B_{H_f}|-\sum_{x\in
B_{H_f}}d_{G-A_{H_f}}(x)-f(A_{H_f})\\
&=\sum_{i=1}^k( |B_i|+r_i-\sum_{x\in
B_i}d_{G-A_{H_f}}(x))-f(A_{H_f})\\
&=\sum_{i=1}^k( |B_i|+r_i-e_G(B_i,V(R_i)-B_i))-f(A_{H_f})\\
&=k-f(A_{H_f})\\
&=o(G[D_{H_f}\cup B_{H_f}])-f(A_{H_f})\\
&\leq o(G-A_{H_f})-f(A_{H_f}),
\end{align*}
a contradiction. This completes the proof. \qed

\noindent \textbf{Proof of Theorem \ref{main_odd}.} By Lemma
\ref{non_consecu}, $G$ is not $H_f$-critical. Then we have
$A_{H_f}\cup B_{H_f}\neq \emptyset$.  For every component $C_i$ of
$G[D_{H_f}]$ and every vertex $v$,  we claim $E_G(C_i,v)\leq 1$.
Otherwise, suppose that  $e_G(C_i,v)\geq 2$. By Lemma
\ref{lem_misses}, we have $d_F(v)\geq 1$ for any $H_f$-optimal graph
$F$, contradicting to $v\in B_{H_f}$.

\medskip
\textbf{ Claim 1.~} $A_{H_f}\neq \emptyset$.
\medskip

  Otherwise, suppose that
$A_{H_f}= \emptyset$. Then we have $B_{H_f}\neq \emptyset$.  Let
$D_{H_f}=C_1\cup\cdots\cup C_k$ and $B_{H_f}=\{v_1,\ldots,v_r\}$.
Since $G$ is connected, then we have $C_{H_f}= \emptyset$ and
$e_G(C_i,B_{H_f})\geq 1$ for $1\leq i\leq k$. Now we show that there
exists a component $C_i$ of $G[D_{H_f}]$ such that
$e_G(C_i,B_{H_f})=1$. Otherwise, assume that $e_G(C_i,B_{H_f})\geq
2$ for $1\leq i\leq k$. Note that $d_G(v)\geq f(v)-1\geq 3$ for all
$v\in B_{H_f}$. Then we have $e_G(D_{H_f},B_{H_f})\geq (2k+3r)/2>
k+r$.
 For any $H_f$-optimal subgraph $F$, by Lemma \ref{lem_misses}, it
misses at most $k$ edges of $E_G(B_{H_f},D_{H_f})$ and so it
contains at least $r+1$ edges of $E_G(B_{H_f},D_{H_f})$. Hence there
exists a vertex $v\in B_{H_f}$, such that $d_F(v)\geq 2$,
contradicting to $v\in B_{H_f}$. Without loss of generality, suppose
that $e_G(C_1,v_1)=1$ and  $u_1v_1\in E(G)$, where $u_1\in V(C_1)$
and $v_1\in B_{H_f}$. By Lemma \ref{lem_D_critical}, $C_1$ is
$H'_f$-critical, where $H'_f:V(C_1)\rightarrow 2^N$ be a set
function such that $H'_f=\{0,2,\ldots,f(v_1)-2,f(v_1)-1\}$ and
$H'_f(u)=H_f(u)$ for all $u\in V(C_1)-v_1$. Note that
$d_{C_1}(v_1)\geq f(v_1)-2$ and $d_{C_1}(u)\geq f(u)-1$ for all
$u\in V(C_1)-v_1$, a contradiction by Lemma \ref{non_consecu}. This
complete Claim 1.

Let $m$ denote the number of components and $G[D_{H_f}\cup
B_{H_f}]=R_1\cup\cdots\cup R_m$.  Suppose that $V(R_i)\cap
B_{H_f}=B_i$ and $R_i$ contains $r_i$ connected components of
$G[D_{H_f}]$. Then by Theorem \ref{lem_main} and Claim 1 of Theorem
\ref{lem_main}, there exists nonempty $A_{H_f}$, such that
 \begin{align*}
0<\nabla(H_f)&=\omega(G[D_{H_f}])+|B_{H_f}|-\sum_{x\in
B_{H_f}}d_{G-A_{H_f}}(x)-f(A_{H_f})\\
&=\sum_{i=1}^m( |B_i|+r_i-\sum_{x\in
B_i}d_{G-A_{H_f}}(x))-f(A_{H_f})\\
&=m-f(A_{H_f})\\
&=o(G[D_{H_f}\cup B_{H_f}])-f(A_{H_f})\\
&\leq o(G-A_{H_f})-f(A_{H_f}),
\end{align*}
a contradiction. This completes the proof. \qed

\end{document}